\documentclass{amsart}
\title[Stable rank and continuous fields]{On the stable rank of 
algebras of operator fields over metric spaces} 

\author{Ping Wong Ng}
\address{Department of Mathematics \\ University of Toronto \\
100 St. George St., Room 4072 \\
Toronto, Ontario\\
M5S 3G3 \\
Canada}
\email{png@math.toronto.edu}
\author{Takahiro Sudo}
\address{Department of Mathematical Sciences \\ Faculty of Science \\
University of the Ryukyus \\
Nishihara-cho, Okinawa \\
903-0213\\
Japan}
\email{sudo@math.u-ryukyu.ac.jp}

\date{\today}

\newtheorem{thm}{Theorem}[section]
 
\newtheorem{lem}[thm]{Lemma}

\theoremstyle{remark}

\newcommand{\B}{\mathcal{B}}
\newcommand{\F}{\mathcal{F}}
\newcommand{\A}{\mathcal{A}}

\newcommand{\T}{\mathbb T}
\newcommand{\C}{\mathbb C}
\newcommand{\M}{\mathcal M}

\begin{document}

\begin{abstract} Let $\Gamma$ be a finitely generated, torsion-free, 
two step nilpotent group.  Let $C^*(\Gamma)$ be the universal $C^*$-algebra
of $\Gamma$.   We show that $acsr( C^*(\Gamma)) = 
acsr( C( (\widehat{ \Gamma })_1)$, where for a unital $C^*$-algebra $\A$,
$acsr(\A)$ is the absolute connected stable rank of $\A$, and where
$(\widehat{\Gamma})_1$ is the space of one-dimensional representations of
$\Gamma$.  For the case of stable rank, we have close results.  In the 
process, we give a stable rank estimate for maximal full algebras of 
operator fields over metric spaces.
\end{abstract}

\maketitle

\section{Introduction}

     Rieffel [\textbf{17}] introduced the notion of 
$\emph{stable rank}$ for $C^*$-algebras
as the noncommutative version of complex dimension of ordinary topological
spaces.  It turns out that the stable rank of a unital $C^*$-algebra is the 
same as its Bass stable rank (see [\textbf{8}]).  

     There has been much work in computing the stable
ranks of the universal $C^*$-algebras of various connected Lie groups.  
The greatest progress has been made 
in the case of type I solvable Lie groups (see [\textbf{20}], [\textbf{22}] 
and [\textbf{23}]).   Roughly speaking,
it has been shown that the stable rank of the universal $C^*$-algebra of 
a type I solvable Lie group $G$ is controlled by the ordinary topological 
dimension of the space of one-dimensional representations of $G$.  

     Recently, the stable ranks of the universal $C^*$-algebras of
a class of non-type I solvable Lie groups (which include the Mautner group)
have been computed (see [\textbf{21}]). 

    In this paper, we compute the stable ranks of the universal $C^*$-algebras
of a class of non-type I amenable discrete groups.  Specifically, our
main result is 

\begin{thm}
Let $\Gamma$ be a finitely generated, torsion-free, two-step nilpotent group.
Let $C^*(\Gamma)$ be the universal $C^*$-algebra of $\Gamma$.  Then
\begin{enumerate}
\item $acsr( C^*( \Gamma )) = acsr(C( (\widehat{\Gamma})_1))$,
\item $sr( C( (\widehat{\Gamma})_1)) \leq sr( C^*( \Gamma )) \leq 
sr( C( (\widehat{\Gamma})_1)) ) + 1$,
\item  if the topological dimension $dim( (\widehat{\Gamma})_1)$ is even,
then $sr(C^*(\Gamma)) = sr( C( (\widehat{\Gamma})_1))$.
\end{enumerate}
\noindent Here $(\widehat{\Gamma})_1$ is the space of one-dimensional   
representations of $\Gamma$.  Also, for a unital $C^*$-algebra $\A$, 
$sr( \A )$ is the stable rank of $\A$, and $acsr(\A)$ is the absolute connected
stable rank of $\A$.  

\end{thm}      \noindent We note that for a unital $C^*$-algebra $\A$, the absolute connected
stable rank of $\A$ is numerically the same as the stable rank of the tensor
product $C[0,1] \otimes \A$.

     A key step in our proof of Theorem 1.1, is the following stable 
rank estimate for algebras of operator fields over metric spaces, which is
of independent interest:

\begin{thm}
Suppose that $X$ is a $\sigma$-compact, locally compact,  
$k$-dimensional metric space.  Suppose that $\A$ is a maximal full algebra of
operator fields over $X$ with fibre algebras, say, $\{ \A_t \}_{ t \in X}$ 
such that $\A_t$ is unital for all $t \in X$.  Then the stable rank of $\A$
satisfies the inequality
\[  
sr( \A ) \leq sup_{ t \in X } sr( C([0,1]^k ) \otimes \A_t ).
\]
\end{thm}
     
    We note that Theorems 1.1 and 1.2 generalize results from [\textbf{13}], 
where
we compute the stable ranks of the universal $C^*$-algebras of the (possibly
higher rank) discrete Heisenberg groups.

     General references for stable rank are [\textbf{14}] and 
[\textbf{17}].
  General references for 
full algebras of operator fields are [\textbf{6}], [\textbf{11}] and
[\textbf{24}].   General references for the
representation theory of finitely generated, two-step nilpotent groups are
[\textbf{2}], [\textbf{15}] and [\textbf{16}] (also see [\textbf{9}]). 

     In what follows, for a $C^*$-algebra $\A$, ``$sr(\A)$'' and ``$acsr(\A)$"
will denote the stable rank and absolute connected stable rank of $\A$ 
respectively.  If, in addition, $\A$ is unital, then for every positive integer
$M>0$, $Lg_M( \A )$ will be the set of all $M$-tuples $(a_1, a_2, ..., a_M)$ in
$\A^M$ such that $\sum_{j=1}^M (a_j)^* a_j$ is an invertible element of 
$\A$.
     Also, for a metric space $X$, for a point $x \in X$ and real number
$r > 0$, ``$B(x,r)$" will denote the open ball of radius r (with respect to the
metric on $X$) about $x$. 

\section{Main results}
     In [\textbf{13}], we proved the following result:
\begin{thm}  Suppose that $\A$ is a unital maximal full algebra of operator 
fields with base space the $k$-cube $[0,1]^k$ and fibre algebras, say,
$\{ \A_t \}_{ t \in [0,1]^k}$.  Then the stable rank of $\A$ satisfies the
inequality
\[
sr( \A ) \leq sup_{t \in [0,1]^k } sr( C([0,1]^k) \otimes \A_t).
\]
\end{thm}
The key technique within the proof of Theorem 2.1, was the following technical
result, which 
we state as a lemma:
\begin{lem}
Suppose that $\A$ is a unital maximal full algebra of operator fields with
base space $[0,1]$ and fibre algebras, say, $\{ \A_t \}_{t \in [0,1]}$.  
Suppose that   $M =_{df} sup_{t \in [0,1]} sr( C[0,1] \otimes \A_t )$ is 
a finite number.
Let $q$, $r$ be real numbers, with $0 < q < r < 1$, and let 
$\A([0,r])$ and $\A([q,1])$ be the restrictions of the operator fields in
$\A$ to $[0,r]$ and $[q,1]$ respectively.  Now  let $\epsilon>0$ be given
and suppose that for $j=1,2,...,M$, 
$\{ a_j(t) \}_{ t \in [0,r]}$ is an operator field in $\A([0,r])$ and 
$\{ b_j(t) \}_{ t \in [q,1]}$ is an operator field in $\A([q,1])$ such that 
\begin{enumerate}
\item[(a)] $\| a_j(t) - b_j(t) \| < \epsilon$ for all $t \in [q,r]$ and for
$j=1,2,..,M$,
\item[(b)] $\sum_{j=1}^M (a_j(t))^* a_j(t)$ is an invertible element of
$\A_t$ for all $t \in [0,r]$, and
\item[(c)] $\sum_{j=1}^M (b_j(t))^* b_j(t)$ is an invertible element of 
$\A_t$ for all $t \in [q,1]$.
\end{enumerate}
\noindent Then there are operator fields $\{ c_j(t) \}_{t \in [0,1]}$ in 
$\A$, $j=1,2,...,M$,  such that  
\begin{enumerate}
\item $\| c_j(t) -a_j(t) \| < \epsilon$ for all $t \in [0,r]$ and for    
$j=1,2,..., M$,
\item $\| c_j(t) - b_j(t) \| < \epsilon$ for all $t \in [q,1]$ and for 
$j=1,2,...,M$, and 
\item $\sum_{j=1}^M (c_j(t))^* c_j(t)$ is invertible in $\A_t$ for all
$t \in [0,1]$.
\end{enumerate}
\end{lem}
 
\begin{proof}[Proof of Theorem 1.2] 
     By [\textbf{12}] IV.7 page 85 last paragraph and [\textbf{19}] Theorem
1.1, we may assume that 
we have a metric on $X$ such that 
for every point $x \in X$ and for every real number $r>0$, the boundary of
the open ball $B(x,r)$ (with respect to this new metric) is at most 
$k-1$-dimensional.   Henceforth, we will 
be working with this metric. 

     Suppose that $X$ is noncompact. Then  let $X_{\infty}$ be the one-point
compactification of $X$, with point at infinity $\infty$.  We may view $\A$
as a maximal full algebra of operator fields with base space $X_{\infty}$ 
and fibre algebras $\{ \A_t \}_{t \in X_{\infty} }$, where 
$\A_t$ is the
same as before when $t \neq \infty$, and $\A_{\infty} = \{ 0 \}$ the zero
$C^*$-algebra.  
Let $\A^{+}$ be the unitization of $\A$. 
By [\textbf{10}] Theorem 1 and Corollary 1,  $\A^{+}$ is a unital maximal 
full algebra of
operator fields with base space $X_{\infty}$ and fibre algebras 
$\{ (\A_t)_1 \}_{t \in X_{\infty}}$ where $(\A_t )_1 = \A_t$ for 
$t \neq \infty$ and $(\A_{\infty})_1 = \mathbb{C}$ (the complex numbers).
A continuity structure 
$\widetilde{\F}$ is the set of all operator fields of the form
$a + \alpha 1$, where $a$ is in $\A$ and $\alpha$ is a complex number. 
If $X$ is compact, then $\A$ will automatically be unital, and (in the 
arguments that follow) we let
$X_{\infty} = X$ and $\A^{+} = \A$, and we need not consider the point
$\infty$ at all. 

     Now suppose that $M=_{df} sup_{t \in X} sr( C([0,1]^k) \otimes \A_t )$ 
is a finite number.     
Let $(a_1, a_2, ..., a_M)$ be an $M$-tuple in $(\A^+)^M$ and let 
$\epsilon > 0$ be given.   By adding a small scalar multiple of the 
unit if necessary,
we may assume that $(a_1, a_2, ..., a_M)$ is nonzero at $\infty$ 
(noncompact) case.  We may also assume that $\epsilon$ is small enough so that 
for any other $M$-tuple $(c_1, c_2, ..., c_M)$, if 
$c_j (\infty)$ is within $\epsilon$ of $a_j (\infty)$ for all $j$, 
then the $M$-tuple
$(c_1, c_2, ..., c_M)$ is also not the zero vector at $\infty$.

     Now we can choose a sequence of nonempty open balls 
$\{ B(x_i, r_i) \}_{i=1}^{\infty}$ in $X$ and a sequence of $M$-tuples
$\{ (f_{i,1}, f_{i,2},..., f_{i,M} ) \}_{i = 1}^{\infty}$ in $(\A^+)^M$ such
that 
\begin{enumerate}
\item[(a)] $X$ is covered by the union of all the open balls 
$B(x_i, r_i)$, $r = 1,2
......$,  
\item[(b)] for every $i$, $i=1,2,3...$, there is a strictly positive number
$\delta_i$ such that $\sum_{j=1}^M f_{i,j}(t)^* f_{i,j}(t)$ is an invertible
element of $\A_t$ for all $t \in B( x_i, r_i + \delta_i)$, and  
\item[(c)] there is an increasing sequence of integers
$\{ N_n \}_{n = 1}^{\infty}$ such that $f_{i,j}$ is within $\epsilon/{2^n}$
 of $a_j$ for $i \geq N_n$, $i=1,2,3,....$ and $1 \leq j 
\leq M$.  
\end{enumerate}
\noindent  Condition (a) uses the $\sigma$-compactness of $X$.
 Condition (b) requires the use of existence and continuity of 
operator fields in full algebras of operator fields (see the definition of
full algebras of operator fields in [\textbf{6}], [\textbf{10}] or 
[\textbf{24}]).  One also needs to use the fact  
that in a unital $C^*$-algebra, any element close enough to the unit 
is invertible.
Condition (c) requires the maximality of the full algebra of operator fields
(see [\textbf{10}] Proposition 1 and [\textbf{24}] Theorem 1.1) as well as the
$\sigma$-compactness of $X$. 
\noindent Henceforth, we let ``(+)" denote properties (a) - (c) collectively.

     By the $\sigma$-compactness of $X$, we may additionally assume that 
such that for every $n$, if $i \leq N_n$ and $j > N_n$ then 
$r_j < (1/8)r_i$, for all $i,j$, i.e., the size of the open balls $B(x_i, r_i)$
are, approximately, ``decreasing uniformly" with rate 
$1/8$.  We may also assume that for each $i$, $i=1,2,3,...$, 
the closure of $B(x_i, r_i + \delta_i)$ is a compact subset of $X$. 

      Our procedure for constructing an $M$-tuple in $Lg_M( \A^+ )$ which
will approximate $(a_1, a_2, ..., a_M)$ to within $\epsilon$ is to construct
$M$ sequences of operator fields $\{ \alpha_{j}^n \}_{n=1}^{\infty}$ 
$j=1,2,..., M$, which satisfy the following conditions:
\begin{enumerate}
\item $\alpha_{j}^n$ is an operator field over 
$\bigcup_{i=1}^{N_n} \overline{B(x_i, r_i)}$ for 
$n=1,2,3,...$ and for $1 \leq j \leq M$,
\item $\alpha_{j}^n$ is within $\epsilon/2$ of $a_j$ 
(over $\bigcup_{i=1}^{N_n} \overline{B(x_i, r_i)}$), 
for $n=1,2,3,...$ and
for $1 \leq j \leq M$,
\item $\alpha_{j}^n$ is within $\epsilon/{2^n}$ of $a_j$ over  
$\overline{B(x_{N_n}, r_{N_n})}$, 
\item $\sum_{j=1}^M (\alpha_{j}^n (t) )^* \alpha_{j}^n (t)$ is invertible
in $\A_t$ for all $t \in \overline{B(x_i, r_i)}$ and for $i \leq N_n$, and
\item for $i \leq m \leq n$, $(\alpha_{1}^m, \alpha_{2}^m,..., \alpha_{M}^m) 
= (\alpha_{1}^n, \alpha_{2}^n, ..., \alpha_{M}^n)$ over the ball
$B(x_i, r_i / 2)$.
\end{enumerate}
\noindent We let ``(*)" denote conditions (1) - (4) collectively.

     For simplicity, let us assume that for every integer $n$, $N_n = n$.     
We now construct the operator fields $\{ \alpha_{j}^n \}_{n=1}^{\infty}$,
$1 \leq j \leq M$, recursively on $n$ (for all $j$ at each step $n$). 
For $n=1$, just let  
$( \alpha_{1}^1, \alpha_{2}^1, ..., \alpha_{M}^1 ) = ( f_{1,1}, f_{1,2}, ...,
f_{1,M})$.
Now suppose that $(\alpha_{1}^{n}, \alpha_{2}^{n}, ..., \alpha_{M}^n)$ has
been constructed.  To construct
$(\alpha_{1}^{n+1}, \alpha_{2}^{n+1}, ..., \alpha_{M}^{n+1})$, we need to 
``connect" $(f_{n+1, 1}, f_{n+1, 2},..., f_{n+1, M})$ with  
$(\alpha_{1}^{n}, \alpha_{2}^{n}, ..., \alpha_{M}^n)$  over an appropriate 
subset of $X$.    
We may assume 
that $\bigcup_{i=1}^{n} \overline{B(x_i, r_i)}$
is nonempty 
(for otherwise, it would be immediate).   

    Let $d$ be the positive real number which is the minimum of the 
quantities $\delta_{n+1}$ and $r_{n+1}$. 
Let $F$ be the set of all points $x$ in $\bigcup_{i=1}^n 
\overline{B(x_i, r_i)}$ whose distance from $x_{n+1}$ is between (and 
including) 
$r_{n+1}$ and $r_{n+1} + d$.    
For $s \in [0,1]$, let $F_s$ be the set of points in
$F$ which have distance $(1-s)r_{n+1} + s(r_{n+1} + d)$ from $x_{n+1}$.  
Let
$\A(F)$ be the $C^*$-algebra gotten by
taking the restriction of (the operator fields in) $\A$ to $F$.  Then   
 $\A(F)$ can be realized as  a unital maximal full algebra of operator 
fields with base space $[0,1]$ and fibre algebras, say, 
$\{ B_s \}_{s \in [0,1]}$.  For each $s \in [0,1]$, the fibre algebra 
$\B_s$ is the restriction of $\A$ to $F_s$, and for each element $a \in \A(F)$,
its fibre at $s \in [0,1]$
(with respect to this continuous field representation) is the 
restriction of $a$ to $F_s$.  Continuity and maximality follows from the 
continuity and maximality of the algebra of operator fields $\A$.  

     Therefore, $C[0,1] \otimes \A(F)$ can be realized as a unital maximal
full algebra of operator fields with base space $[0,1]$ and fibre algebras
$\{ C[0,1] \otimes \B_s \}_{s \in [0,1]}$.  The continuity structure consists
of all operator fields of the form 
$s \mapsto \sum_{i=1}^N f_i \otimes b_i (s)$, where the $f_i$s are in 
$C[0,1]$ and the $b_i$s are continuous operator fields in $\A(F)$ (with 
respect to the continuous field decomposition of $\A(F)$ in the previous 
paragraph).  Hence, by Theorem 2.1, the stable rank of $C[0,1] \otimes \A(F)$   
satisfies $sr( C[0,1] \otimes \A(F)) \leq sup_{s \in [0,1]} sr( C[0,1] \otimes
\B_s )$. 
     
     But for $s \in [0,1]$, $\B_s$ can be realized 
as a unital maximal full
algebra of operator fields with base space $F_s$ (a compact
metric space) and 
fibre algebras $\{ \A_t \}_{ t \in F_s}$ (since $\B_s$ is the restriction of
$\A$ to $F_s$).   Hence, for $s \in [0,1]$, $C[0,1] \otimes \B_s$ can be 
realized as a unital maximal full algebra of operator fields with base 
space $F_s$ and fibre algebras $\{ C[0,1] \otimes \A_t \}_{t \in F_s }$.
But by our assumption on the metric in the first paragraph of this
proof, $F_s$ is a metric space with dimension less than or equal to $k-1$. 
Hence,
 we have, by induction,
that $sr( C[0,1] \otimes \B_s ) \leq sup_{t \in F_s} sr( C([0,1]^{k-1} 
\otimes C[0,1] \otimes \A_t )$ (The induction is on the dimension of the
base space.  Note that when $F_s$ is zero-dimensional,  the 
stable rank estimate will be immediate, since we can choose a 
\emph{finite, clopen} covering for $F_s$, which satisfies the properties 
in (+).
Hence the base case is immediate).  
From this and the previous paragraph,
$sr( C[0,1] \otimes \A(F) ) \leq M$.

        By Lemma 2.2, it follows that there is an $M$-tuple of operator fields
$(\alpha_{1}^{n+1}, \alpha_{2}^{n+1}, ..., \alpha_{M}^{n+1})$  on 
$\overline{B(x_{n+1} , r_{n+1})} \cup \bigcup_{i=1}^n \overline{B(x_i, r_i)}$
such that  
$\alpha_{j}^{n+1} = f_{n+1, j}$ on $B(x_{n+1}, r_{n+1})$,  
$\alpha_{j}^{n+1} = \alpha_{j}^n$ on $\bigcup_{i=1}^M B(x_i, r_i) -
F$, and the $\alpha_{j}^{n+1}$s satisfy (1), (2) and (4) in (*).  
Condition (5) in (*) is satisfied, since $d$ was chosen to be less than
or equal to $r_{n+1}$, and the latter is strictly less than $(1/8)r_n$.
Finally, condition (3) in (*) is satisfied since $f_{n+1,j}$ is within
$\epsilon/{2^{n+1}}$ of $a_j$.
 
     In the general case where $N_{n+1}$ is not necessarily equal to $n+1$, we
need to repeat the preceding procedure a finite number of times, 
in the natural way, in order to go from $(\alpha_{1}^n, \alpha_{2}^n, ..., 
\alpha_{M}^n)$ to $(\alpha_{1}^{n+1}, 
\alpha_{2}^{n+1}, ..., \alpha_{M}^{n+1})$. 

      Also, when $X$ is compact, the preceding procedure will stop at 
finitely many steps and the sequences $\{ \alpha_{j}^n \}_{n =1}^{\infty}$,
$1 \leq j \leq M$, will all be finite.  We leave to the reader the obvious 
modifications that need to be made.
  
    Now suppose that we have constructed sequences of operator fields
$\{ \alpha_{j}^n \}_{n=1}^{\infty}$, $j = 1, 2, ..., M$ as in (*).
We then construct an $M$-tuple of continuous operator fields in 
$(\A^+)^M$ as follows:
let $\alpha_j (t) = \alpha_{j}^n (t)$ for $t \in B(x_n, {r_n}/2)$ and (in the
noncompact case) let $\alpha_j (\infty) = a_j (\infty)$.  Continuity at the 
point $\infty$ is ensured by condition (3) in (*).   
Then $\alpha_j$ is within $\epsilon /2$ of $a_j$ for $j=1,2,...,M$.
Moreover, that 
$ \sum_{i=1}^M (\alpha_j(\infty))^* \alpha_j (\infty)$ is invertible follows
from the invertibility of 
$\sum_{i=1}^M (a_j (\infty))^* a_j (\infty)$ and the smallness of the 
the $\epsilon$, both of which were assumed at the beginning.
Hence,  $(\alpha_1, \alpha_2, ..., \alpha_M) \in Lg_M ( \A^+)$ and 
$\alpha_j$ is within $\epsilon$ of $a_j$ for $j=1,2,..., M$. 
 
 \end{proof} 

     The result of the next computation is surely known (See [\textbf{18}],
the comments after the proof of Proposition 3.10).
\begin{lem}
If  $\A_{\Theta}$ is a simple noncommutative torus and $\T^k$ the ordinary
$k$-torus,
then $sr( C(\T^k) \otimes \A_{\Theta}) = 2$.
\end{lem}
\begin{proof}
By [\textbf{3}] Theorem 1.5 and [\textbf{17}] the proof of Corollary 7.2, 
$sr( C(\T^k) \otimes \A_{\Theta})$ is a finite number.  Hence by
[\textbf{17}] Theorem 6.1, let $l$ be a positive integer such that both
$sr( \mathbb{M}_{2^l}(\C) \otimes C(\T^k) \otimes \A_{\Theta})$ and 
$sr( \mathbb{M}_{3^l}(\C) \otimes C(\T^k) \otimes \A_{\Theta})$
are less than or equal to $2$.   Let $\A_{\Theta} = 
\overline{\bigcup_{n=1}^{\infty} \A_n}$ be the inductive limit decomposition
of $\A_{\Theta}$ given in [\textbf{3}] Corollary 2.10.  

     Now let a positive real number $\epsilon > 0$ and a positive integer
$m > 0$ be given.  Let $a_1$ and $a_2$ be arbitrary elements of $C(\T^k) 
\otimes \A_m$.  Choose an integer $n > m$ such that there are $(b_1, b_2) 
\in Lg_2 ( \mathbb{M}_{2^l}(\C) \otimes C(\T^k) \otimes \A_n)$ and       
$(c_1, c_2) \in Lg_2 (  \mathbb{M}_{3^l}(\C) \otimes C(\T^k) 
\otimes \A_n)$, with $b_j$ within $\epsilon$ of $a_j \otimes 1_{\M_{2^l}}$ 
and $c_j$ within  
$\epsilon$ of $a_j \otimes 1_{\M_{3^l}}$, $j = 1,2$. 
Then it follows from the proof of [\textbf{3}] Corollary 2.10, that we 
can choose 
an integer $N > n$ and choose a  
finite dimensional subalgebra $\B \subseteq \A_N$ such that there exists
$(d_1, d_2) \in Lg_2 (C(\T^k) \otimes C^*(\A_n, \B))$ with $d_j$ being within
$\epsilon$ of $a_j$ for $j=1,2$.  
But $\epsilon$, $m$ and $a_j$ were arbitrary.  Hence, $sr( C(\T^k) \otimes 
\A_{\Theta}) \leq 2$. 

     Now $K_1 (\A_{\Theta}) = \mathbb{Z}^{2^{p-1}} \neq 0$ where $p$ is the
dimension of the noncommutative torus $\A_{\Theta}$ (i.e., $\A_{\Theta}$ is
a noncommutative $p$-torus).  So we can find a positive integer $n$ such
that $GL_n(\A_{\Theta}) \neq GL_n (\A_{\Theta})_0$, where
$GL_n(\A_{\Theta})$ is the group of invertibles in $\M_n(\A_{\Theta})$ and
$GL_n (\A_{\Theta})_0$ is the connected component of the identity in
$GL_n(\A_{\Theta})$.  Therefore, the connected stable rank
$csr(\M_n(\A_{\Theta}) \geq 2$. But
$csr(\M_n(\A_{\Theta}) \leq 
sr(\mathbb{M}_{n}(\C) \otimes C(\T^k) \otimes \A_{\Theta}) 
\leq sr(  C(\T^k) \otimes \A_{\Theta})$.  Hence,
$sr( C(\T^k) \otimes \A_{\Theta}) \geq 2$.
\end{proof}

\begin{proof}[Proof of Theorem 1.1]  
 Since $C((\widehat{\Gamma})_1)$ is naturally a quotient of $C^*(\Gamma)$, we
must have that $sr(C^*(\Gamma)) \geq sr(((\widehat{\Gamma})_1))$ and
$acsr(C^*(\Gamma)) \geq acsr(((\widehat{\Gamma})_1))$.  

     By [\textbf{2}] page 390 last paragraph and page 391 first paragraph,
and by [\textbf{15}] Theorem 1.2,  
$C^*(\Gamma)$ can be realized as a unital maximal full algebra
of operator fields with base space $\widehat{Z(\Gamma)}$ and fibre algebras,
say, $\{ \A_{\lambda} \}_{\lambda \in \widehat{ Z(\Gamma )} }$.  Here,
$Z(\Gamma)$ is the centre of $\Gamma$, and 
$\widehat{ Z(\Gamma)}$ is the Pontryagin dual of the centre of $\Gamma$.  
Moreover,  the continuous open surjection corresponding to this 
continuous field decomposition of $C^*(\Gamma)$ is the map $p: Prim(C^* (
\Gamma )) \rightarrow \widehat{ Z(\Gamma)}$ which brings a primitive ideal
of $C^*(\Gamma)$ to its restriction to $Z(\Gamma)$.  

     Also, by [\textbf{2}] page 390 last paragraph and page 391 first
paragraph, by [\textbf{15}] Theorem 1.2, and by [\textbf{16}] Theorem 1, 
for fixed $\lambda \in \widehat{Z(\Gamma)}$,
$\A_{\lambda}$ (as in the previous   
paragraph) can in turn be realized as a unital maximal full algebra of 
operator fields with base space of the form $\T^g \times T$ for some 
commutative $g$-torus $\T^g$ and finite set $T$.  The integer $g$ is less 
than or equal to the rank of $\Gamma/ Z(\Gamma)$.  The fibre algebras are
all isomorphic.  Let $\B_{\lambda}$ be the unique $C^*$-algebra which
all the fibre algebras  are isomorphic to.  Then 
$\B_{\lambda}$ will be either of the
form $\M_n( \C )$ (a full matrix algebra) or $\M_n (\C ) \otimes \A_{\Theta}$
where $\A_{\Theta}$ is a simple noncommutative torus (the former case will 
occur if $p^{-1}(\lambda)$ consists of $n$-dimensional representations and 
the latter will occur 
if $p^{-1}(\lambda)$ consists of infinite dimensional representations).  
We note that 
$g$, $T$ and $\B_{\lambda}$ will all depend on $\lambda$.

 Now let $\Gamma^{(2)}$ be the commutator subgroup of $\Gamma$ (i.e., the 
subgroup of $\Gamma$ generated by elements of the form $xy x^{-1} y^{-1}$ where
$x,y \in \Gamma$).  Let $\Gamma^{(2)}_s$ be the saturation of 
$\Gamma^{(2)}$ (i.e., the smallest sugroup $H$ of $\Gamma$ containing
$\Gamma^{(2)}$ such that for every $x \in \Gamma$, if $x^n \in H$ for some 
strictly positive integer $n$ then $x \in H$). 
Since $Z(\Gamma)$ (the centre of $\Gamma$) is a saturated 
subgroup of $\Gamma$ (i.e., $x \in \Gamma$ and $x^n \in Z(\Gamma)$ for some 
strictly positive integer $n$ implies that $x \in Z(\Gamma)$), 
$\Gamma^{(2)}_s$ is a saturated subgroup of $Z(\Gamma)$.
Hence, we have a decomposition $Z(\Gamma) = \Gamma^{(2)}_s \oplus F$, where
$F$ is a saturated free abelian subgroup of $Z(\Gamma)$.  This in turn gives
a decomposition $\widehat{Z(\Gamma)} = \widehat{\Gamma^{(2)}_s} \times
\widehat{F}$.  

     Now let $N$ be a positive integer such that for all $n \geq N$, 
$sr( \M_n (\C) \otimes C( \T^{h+1} ) ) \leq 2$, where $h$ is the rank of 
$\Gamma$.  Let $S$ be the set of  all 
$\lambda \in \widehat{Z(\Gamma)}$ 
such that $p^{-1}(\lambda)$ consists of $m$-dimensional representations with
$m \leq N$.  With respect to the decomposition of $\widehat{Z(\Gamma)}$ given
in the previous paragraph, $S$ must have the form $\{ \lambda_1, \lambda_2,...,
\lambda_k \} \times \widehat{F}$, for a finite set of points $\lambda_i \in 
\widehat{\Gamma^{(2)}_s}$.  (Suppose that $x_1, x_2, ..., x_q$ are elements of
$\Gamma$ so that $x_1 / Z(\Gamma), x_1 / Z(\Gamma), ..., x_q / Z(\Gamma)$ give
a basis for $\Gamma /Z(\Gamma)$.  Suppose that $\pi$ is an $m$-dimensional 
representation of $\Gamma$.  Then 
the scalar values 
$\pi (x_i x_j (x_i)^{-1} (x_j)^{-1})$,  for $ 1 \leq i, j \leq q$, must all
be rational numbers  which can be placed in the form $r/q$ where $q \leq m$.
These scalar values determine the values of $\pi$ on $\Gamma^{(2)}$ and there
are only finitely many possibilities for them.  And since $\Gamma^{(2)}$ has
finite index in $\Gamma^{(2)}_s$, there are only finitely many possibilities
for the restriction of $\pi$ to $\Gamma^{(2)}_s$.)

     Let $J$ be the ideal of $C^*(\Gamma)$ consisting of all operator fields
which vanish on $S$.  Then $J$ is a maximal full algebra of operator fields
with base space $\widehat{Z(\Gamma)} - S$.  The quotient $C^*(\Gamma) / J$ is
a unital maximal full algebra of operator fields with base space $S$. Indeed,
$C^*(\Gamma) / J$ is the restriction, to $S$, of the operator fields in 
$C^*(\Gamma)$.  
Now from the exact sequence
$
0 \rightarrow C[0,1] \otimes J \rightarrow 
C[0,1] \otimes C^*(\Gamma) \rightarrow C[0,1] \otimes 
C^*(\Gamma)/ J \rightarrow
0,
$
\noindent and by [\textbf{7}] Corollary 2.22, 
we get that $sr(C[0,1] \otimes C^*(\Gamma) ) = max \{ 
sr( C[0,1] \otimes J ) , sr ( C[0,1] \otimes C^*(\Gamma)/J ) \}$.  By 
Theorem 1.2, 
and our definitions of $N$ and $S$, $sr( C[0,1] \otimes J ) \leq 
sup_{ \lambda \in \widehat{Z(\Gamma)} - S} sr( C[0,1] \otimes \A_{\lambda})$.
But for $\lambda \in \widehat{Z(\Gamma)} - S$, we have, by the definitions of
$N$ and $S$, by Lemma 2.3, by [\textbf{17}] Proposition 1.7 and Theorem 
6.1, and by our discussion of the 
continuous field decomposition
of $\A_{\lambda}$, that $sr( C[0,1] \otimes \A_{\lambda} ) \leq 2$.  Hence,
$sr( C[0,1] \otimes J ) \leq 2$.  

     Now by Theorem 1.2, $sr( C[0,1] \otimes C^*(\Gamma) / J ) \leq 
sup_{\lambda \in S } sr( C[0,1] \otimes \A_{\lambda} )$.  
But for $\lambda \in S$, we have that the continuous field decomposition of
$\A_{\lambda}$ has base space $\T^g \times T$ where $g$ is less than the rank
of $\Gamma / Z(\Gamma)$ and $T$ is a finite set.  The fibre algebras are 
all isomorphic to $\M_m(\C)$, for some integer $m$.  Hence, by Theorem 1.2,
for $\lambda \in S, sr( C[0,1] \otimes C^*(\Gamma))  \leq  
sr( C(\T^{g + 1} ) )  = sr( C[0,1] \otimes C( (\widehat{ \Gamma })_1 )$.  
From this and the previous paragraph, 
$acsr( C^*( \Gamma ) ) = acsr( C( ( \widehat{ \Gamma })_1 ))$.

     The proofs of the other statements of the theorem now follow from our
computation for absoluted connected stable rank.  

     By [\textbf{17}] Theorem 4.3
 and our result for absolute connected stable rank, we 
have that  $sr( C^*( \Gamma ) ) \leq 
sr( C[0,1] \otimes C^*(\Gamma))
= sr( C[0,1] \otimes C((\widehat{\Gamma})_1)))$ 

     But by [\textbf{17}] Corollary 7.2,   $sr( C[0,1] \otimes C^*(\Gamma) ) 
\leq sr( C^*(\Gamma)) + 1$.  Hence, $sr( C( (\widehat{\Gamma})_1))) \leq
sr( C^*( \Gamma )) \leq sr( C( (\widehat{\Gamma})_1))) + 1$.  

    Also, by [\textbf{17}] Proposition 1.7, if 
$dim( (\widehat{\Gamma})_1 )$ is even, then 
$sr( C[0,1] \otimes C( (\widehat{\Gamma})_1) ) = sr(C((\widehat{\Gamma})_1))$. 
Hence, if 
$dim( (\widehat{\Gamma})_1 )$ is even, then $sr(C^*(\Gamma)) = 
sr( C((\widehat{\Gamma})_1))$.   
\end{proof}
       
\section*{References} 

\newcounter{bean}
\begin{list}{\textbf{\arabic{bean}}}{\usecounter{bean}}             

\item   \textsc{J. Anderson and W. Paschke}, 
`The rotation algebra',
\emph{Houston J. Math.}, 15, (1989), 1 - 26. 

\item \textsc{L. Baggett and J. Packer}, 
`The primitive ideal space of two-step nilpotent group $C^*$-algebras',
\emph{J. Funct. Anal.}, 124, (1994), 389-426.

\item \textsc{B. Blackadar and  A. Kumjian and M. Rordam}, 
`Approximately central matrix units and the structure of noncommutative tori',
K-Theory, 6, (1992), 267-284.
 
\item \textsc{ K. R. Davidson},
\emph{$C^*$-algebras by example},
Fields Institute Monographs, 6.

\item \textsc{G. A. Elliott and D. E. Evans}, `The structure of
the irrational rotation $C^*$-algebra', \emph{Ann. of Math. (2)}, 138,
no. 3, (1993), 477-501.

\item \textsc{J. M. G. Fell},
`The structure of algebras of operator fields', \emph{Acta Math.}, 
106, (1961), 233-280.

\item \textsc{N. E. Hassan}, `Rangs Stables de certaines 
extensions', \emph{J. London Math. Soc.}, (2), 52, (1995), 605-624.

\item \textsc{R. Herman, L. N. Vaserstein}, 
`The stable rank of $C^*$-algebras', \emph{Invent. Math.}, 
77, (1984), 553-555.

\item \textsc{R. Howe}, `On representations of discrete, finitely generated,
torsion-free, nilpotent groups', \emph{Pacific J. Math}., 
73, (1977), 281-305.

\item \textsc{R-Y Lee}, `On $C^*$-algebras of operator fields',
\emph{Indiana U. Math. Journal}, 25, no. 4, (1976), 303-314.  

\item  \textsc{S. T. Lee and J. Packer},
`Twisted group $C^*$-algebras for two-step nilpotent and generalized
discrete Heisenberg group', \emph{J. Operator Th.}, 34, no. 1, (1995), 
91-124.

\item \textsc{J. Nagata}, \emph{Modern dimension theory},
revised edition, Sigma Series in Pure Mathematics Volume 2, 
(Heldermann Verlag, Berlin 1983).   

\item \textsc{P. W. Ng and T. Sudo} 
`On the stable rank of algebras of operator fields over an $N$-cube'.
preprint.  E-print available on the Web at 
http://arXiv.org/abs/math.OA/0203115. 

\item \textsc{V. Nistor},
`Stable range for the tensor products of extensions of $K$ by 
$C(X)$', \emph{J. Operator Th.}, 16, (1986), 387-396.

\item \textsc{J. Packer and I. Raeburn}, 
`On the structure of twisted group $C^*$-algebras',
\emph{Trans. AMS}, 334, no. 2, (1992), 685-718.

\item \textsc{D. Poguntke}, 'Simple quotients of group $C^*$-algebras for 
two-step nilpotent groups and connected Lie groups', Ann. Sci. Ec. Norm. 
Sup., 16, (1983), 151-172.

\item \textsc{M. A. Rieffel},
`Dimension and stable rank in the $K$-theory of $C^*$-algebras',
\emph{Proc. London Math. Soc.}, (3), 46, (1983), 301-333.

\item \textsc{M. A. Rieffel},
`Projective modules over higher-dimensional noncommutative tori',
\emph{Can. J. Math.}, 40, no. 2, (1988), 257-338.
 
\item \textsc{J. J. Roberts}, `A theorem on dimension', 
\emph{Duke Math. J.}, 8, (1941), 565-574.

\item \textsc{T. Sudo}, `Dimension theory of group $C^*$-algebras of 
connected Lie groups of type I',  \emph{J. Math. Soc. Japan}, 52, no. 3,
(2000), 583-590.    

\item \textsc{T. Sudo}, 'Structure of group $C^*$-algebras of Lie semi-direct
products $\mathbb{C} \times_{\alpha}  \mathbb{R}$', 
\emph{J. Operator Theory}, 46, (2001), 25-38.

\item \textsc{T. Sudo and H. Takai}, `Stable rank of the 
$C^*$-algebras of nilpotent Lie groups', \emph{Internat. J. Math.}, 6,
no. 3, (1995), 439-446.

\item \textsc{T. Sudo and H. Takai}, `Stable rank of the 
$C^*$-algebras of solvable Lie groups of type I', \emph{J. Operator Theory},
38, no. 1, (1997), 67-86.

\item \textsc{J. Tomiyama}, 
`Topological representation of $C^*$-algebras',
\emph{Tohoku Math. J.}, (2), 14, (1962), 187-204.
\end{list}

\end{document}